\documentclass{ifacconf}

\usepackage{graphicx}      
\usepackage{natbib}        
\usepackage{amsmath}
\usepackage{color}
\usepackage{amssymb}  
\usepackage{dsfont}
\usepackage{ mathabx}
\DeclareTextSymbol{\deg}{T1}{6}
\DeclareTextSymbol{\deg}{OT1}{23}

\begin{document}
\begin{frontmatter}
\title{Robust MPC for temperature management on electrical transmission lines}

\thanks[footnoteinfo]{This work has been supported by the RTE-CentraleSup{\'e}lec Chair "The Digital Transformations of Electrical Networks"(https://rtechair.fr/)}

\author[First]{Cl{\'e}mentine Straub} 
\author[First]{Sorin Olaru} 
\author[Second]{Jean Maeght}
\author[Second]{Patrick Panciatici}

\address[First]{Laboratory of Signals and Systems (L2S), CentraleSup{\'e}lec, 
   Universit{\'e} Paris-Saclay, France (e-mail: firstname.lastname@centralesupelec.fr).}
\address[Second]{French transmission system operator, R{\'e}seau de Transport d'Electricit{\'e} (RTE) (e-mail: firstname.lastname@rte-france.com)}

\begin{abstract} 
In the current context of high integration of renewable energies, maximizing infrastructures capabilities for electricity transmission is a general need for Transmission System Operators (TSO). The French TSO, RTE, is developing levers to control power flows in real-time: renewable production curtailment is already employed and large battery storage systems are planned to be installed for congestion management in early 2020. The combination of these levers with the use of Dynamic Line Rating (DLR) helps exploiting the lines at the closest of their limit by managing their temperature in real-time. Unnecessary margins can be reduced, avoiding congestion and excessive generation curtailment. In particular, there is a possible interesting correlation between the transits increase due to high wind farms generation and the cooling effect of wind on power lines in the same area. In order to optimize the electrical transmission network capacities, the present paper advocates the use of a temperature management model, mixing production curtailment and large batteries as control variables. A robust Model Predictive Control framework for local control on electrical lines temperature is presented based on the regulation within tubes of trajectories. Simulations on the French electrical network are conducted to show the effectiveness of the optimization-based control design.
\end{abstract}


\end{frontmatter}

\section{INTRODUCTION}

As renewable energies become more efficient and cheaper, an expansion of the sector is observed worldwide. The report \cite{BP_RTE} on electricity supply and demand for the French territory foresees a dramatic increase in renewable production, which may lead to congestions. The traditional lever to face this issue was network development. However, due to heavy costs and difficulties in obtaining authorizations, TSOs are looking for new solutions that maximize the utilization of the existing network instead of building new infrastructures.

In general the capacity of a line is defined by its Ampacity, which is the maximal constant current meeting security and safety criteria, maintaining a sufficient clearance and avoiding thermal damage to the conductor. This value has been traditionally determined based on the conductor maximal allowable temperature using a worst-case scenario on weather conditions. However, conductor temperatures, and hence its power transfer capability, are affected by wind and solar conditions. Moreover, there is a correlation between high wind power production and high convective cooling. When the wind is strong, power flows are larger, but the transmission capacity is increased in the same time, see \cite{stephen2012guide}.

Therefore, the use of Dynamic Line Rating which considers dynamic limits evolving according to weather conditions has been proposed in literature, see \cite{foss1990dynamic}. Economic DLR gains are highlighted in several papers, as \cite{schell2011using}. Some solutions, giving real-time information on lines ampacities are already available: \emph{Ampacimon} is described  in \cite{cloet2011uprating}. Generation redispatch and operation of large battery systems have also been proposed as an alternative to network reinforcement to deal with congestions in \cite{wen2015enhanced}.

The French TSO defined an operating policy where limitations are evolving seasonally and where overloads are tolerated for short amounts of time as an approximation of the conductor heating, as explained in \cite{mementoRTE}. For example, an overload of 15\% of the worst-case scenario ampacity can be tolerated for 1 minute and an overload of 10\% can be tolerated for 5 minutes. \cite{zonalCongestion2018} designs a real-time controller managing congestion in this operating context. The controller is based on a Model Predictive Control approach combining generation curtailment and large batteries to enforce the limitations.

The aim of the present paper is to go a step further and develop a new management temperature framework acting on the same quick levers (batteries and renewable production curtailment) and monitoring lines temperatures in real-time. The goal of the designed controller is to enforce temperature constraints. The controller thus includes a model of transmission line heating, needs  predictions on temperature evolution and has to handle delays on control. It is natural to design a MPC controller \citep{mayne2000constrained} as operational constraints are included. The method also benefits from the feedback mechanism as linear approximations are considered in the model. It is possible to design a robust controller by means of the Tube MPC strategy \citep{rakovic2012parameterized}, \citep{mayne2005robust} and this is a crucial issue as wind is highly unpredictable.

The main contribution of this paper is to develop a robust MPC controller to manage the real-time transmission lines temperature. A non-linear receding-optimization formulation is considered based on the temperature prediction model. The recursive feasibility will be enforced by widening uncertainties bounds to include the model mismatch and performances will be checked with representative simulations on a test case of the French network.

\section{Modeling} 

\subsection{Electrical transmission network modeling}

The DC load flow principles are used in this article to model power flows evolution \citep{zimmerman2010matpower}. DC load flow consists in a linearization of AC network equations with three assumptions: flat voltage profile, small differences in voltage angles and transmission without losses. The DC modeling enables us to subsequently define  a model-based control problem and exploit the linearity in order to obtain a convex optimization formulation that can be solved very efficiently. The three assumptions made are not major restrictions.The feedback controller will be able to handle them. In order to determine power flows, the knowledge of injections at each node of the network is needed with the DC load flow. Such knowledge on the whole network is an ideal feature but is not necessary: congestion are a local phenomenon and can be handled locally. Power Transfer Distribution Factors (PTDF), described in \cite{bart2005network}, can help representing local power flows, without solving DC load flow on the whole network. The flow representation by  PTDF means and DC load flows are strictly equivalent \citep{zimmerman2010matpower}. Using PTDF, power flows evolution can be described through linear equations. A PTDF is a factor dependent on a line and a bus which gives the amount of power going through the line if the power injection at the bus increases of 1. In other words, let $F_{ij}$ denotes the power flow on line $ij$. Then a change in the flow of line $ij$ $\Delta F_{ij}$ takes place if the injection at bus $n$ increases of $\Delta P_n$, and leads to the formal definition of the factor:
\[ PTDF(ij,n) = \frac{\Delta F_{ij}}{\Delta P_n} \]

The variation $\Delta F_{ij}$ depends on the choice of the 'slack bus'. In a network, the balance between load and generation must be maintained. If the power injection varies at bus $n$, the injection at another bus must vary accordingly and maintain the balance. This second bus is called the slack bus. For a fixed slack bus and a fixed network topology, the rapport $\frac{\Delta F_{ij}}{\Delta P_n}$ is constant. Power flows evolutions between $k$ and $k+1$ time instants can thus be written as:
\begin{equation}
F_{ij}^{k+1} = F_{ij}^{k} + \sum_n PTDF(ij,n)( \Delta u_n^k + w_{n}^k)
\label{power flow}
\end{equation}
where $F_{ij}^k$ is the power flow on line $ij$ at time $k$, $w_{ij}^k$ the  natural change in generation at node $n$ and  
$\Delta u_n^k$ is the injection variation due to control (generation curtailment or battery charge). The sampling time is a model parameter. We consider here a sampling time of approximately 1 minute: as limitations concern temperature and not power flows, 1 minute is precise enough for power flows evolution. 1 minute is also an appropriate sampling time for temperature evolution, as it will be detailed next.

\subsection{Temperature evolution}

Traditional Dynamic Line Rating algorithms give the lines ampacities as a function of weather conditions. The ampacity is updated frequently (every 10 minutes for instance) when new measurements on solar radiation, wind and temperatures (ambient and of the lines) are available. The controller developed here does not have limitations on intensity but only on lines temperatures, as it includes a model of line heating. Equations similar to DLR are used to model this heating. The work in \cite{iglesias2014guide} gives a general method to calculate the thermal rating of overhead lines and presents in particular an algorithm to track conductor temperatures with time considering that weather data and current are provided every 10 minutes. Data are considered constant on the 10 minutes interval and calculations are performed every minute. One minute is small enough compared to the time constants of the line heating and provides accurate result. The algorithm described  in \cite{iglesias2014guide} is adapted to give an equation describing the temperature evolution. The time step used is also one minute. Measurements both on power flows and line heating are available every minute. Equations are based on the discretized heat balance equation for transient state: 
\begin{equation}
    \Delta T = \frac{P_J + P_S - P_C }{m \cdot c} \Delta t
    \label{Delta T}
\end{equation} 
 $m$ is the mass per unit ($kg.m^{-1}$), $c$ is the heat capacity ($J.kg^{-1}.J^{-1}$), $T$ the line temperature, $\Delta t$ the time step, $P_J, P_S$ and $P_C$ respectively the Joule and solar heating, and the convective cooling.
 \begin{equation*}
     P_J = R\cdot I^2
 \end{equation*}
  \begin{equation*}
     P_S = \alpha_S \cdot D \cdot I_T
 \end{equation*}
  \begin{equation*}
     P_C = \pi \cdot \lambda_f \cdot N_u \cdot (T_c - T_a)
 \end{equation*}
$I$ is the total direct current, $R$ the direct current resistance per unit length ($\Omega.m^{-1}$), $\alpha_S\in [0.2,0.9]$ the surface absorbance which depends on the conductor age, $D$ the conductor outer diameter ($m$), $I_T$ the global radiation intensity ($W.m^{-2}$), $\lambda_f$ the air thermal conductivity ($W.K^{-1}.m^{-1}$). $T_c$ / $T_a$ are the conductor / air temperatures and $N_u$ is the dimensionless Nusselt number calculated with the wind speed and direction. The model makes several assumptions: $R$ and $\lambda_f$ are constant with temperature, the conductor temperature is constant across section and the magnetic heating is neglected. The term coming from the Joule heating needs the current intensity I and the DC modeling is based on F, the three phase active power. The relation between I and F is given by: 
\[I = \frac{\sqrt{F^2 + Q^2}}{\sqrt{3} V} \]
with $V$ the phase to phase voltage and $Q$ the reactive power. $Q$ can be measured and $V$ taken to its nominal value. The equation \eqref{Delta T}, replacing $I$ by its expression with $F$, is non-linear due to the square on the active power term $F$.

The evolution of temperature in an electrical line can thus be written as:
\begin{multline}
T^{k+1} =  (1-\frac{\pi \cdot \lambda_f \cdot N_u}{m \cdot c} \Delta t)T^k +\frac{R \cdot \Delta t}{3 V_0^2 \cdot m \cdot c  } (F^k)^2 \\
+ \frac{\Delta t}{m \cdot c}(\alpha_S \cdot D \cdot I_T + \pi \cdot \lambda_f \cdot N_u \cdot T_a + R \frac{Q_0^2}{3 V_0^2})
\label{eqTemp1}
\end{multline}

with $\lambda_f$, $R$, $V_0$, $\alpha_S$, $D$, $m$ and $c$ constants. $N_u$, $Q_0$, $T_a$ and $I_T$ can be regularly updated. Coefficients are time-dependant (variations of weather and electrical values: $N_u$, $Q_0$, $T_a$, $I_T$, $F_0$). However they can be considered constant for prediction purposes as long as their variation is slow with respect to the chosen sampling time.

\section{MPC formulation for temperature management}

\subsection{MPC formulation}

Curtailment actions present time-delay while the battery action is immediate. It is considered that a 1 minute delay is necessary for curtailment (time needed for order transmission). As $\Delta t = 1 min$, the delay for curtailment is one time step. If the delay is larger or the time step smaller, an extended state  space representation can be constructed. The control variables are the power injected in the battery and generation curtailment. Several sites for curtailment are considered, each with a different impact on power flows. This impact is given by the PTDF. The battery is considered strategically situated, at a point where congestion are important. In eq. \eqref{power flow}, power flows are described with variations in generation curtailment and power injected in the battery. The control needs to be expressed in term of differences: $\Delta u_{batt}^k = u_{batt}^{k} - u_{batt}^{k-1}$ for battery injections and $\Delta u_{curt,i}^k = u_{curt,i}^{k} - u_{curt,i}^{k-1}$ for curtailment. $u_{batt}^k \in \mathbb{R}$ is the power injected in the battery at time $k$ and $u_{curt,i}^k \in \mathbb{R}^{n_g}$ the production curtailed in site $i$ at time $k$. $n_g$ is the number of curtailment sites. The control is defined as:
\begin{center}
$u^k=
\begin{pmatrix}
\Delta u_{batt}^k &\quad \Delta u_{curt,i}^k 
\end{pmatrix}$
$\in \mathbb{R}^{n_g+1}$
\end{center}
The state $x$ thus includes $u_{batt}$ and $u_{curt,i}$, as well as the energy stored in the battery $E \in \mathbb{R}$, the power flows on controlled lines $F_{ij} \in \mathbb{R}^{n_l}$ with $n_l$ the number of controlled lines, and their temperature $T_{ij} \in \mathbb{R}^{n_l}$. There are $n$ state variables, with $n = 2(n_l + 1) + n_g$
\begin{center}
$x^k=
\begin{pmatrix}
F_{ij}^k &\quad T_{ij}^k & \quad u_{batt}^k & \quad E_{batt}^k & \quad u_{curt,n}^k
\end{pmatrix} ^T$
$\in \mathbb{R}^{n}$
\end{center}
The dynamical model includes the power flow equations \eqref{power flow}, an integrator to relate $\Delta u$ and $u$, the evolution of the battery state of charge: $E^{k+1}_{batt} = E^{k}_{batt} + \Delta t \cdot u_{batt}^{k}$, and the temperature evolution. A term is added compared to \eqref{eqTemp1} for temperature evolution to take into account the immediate battery action.
\begin{equation}
    T^{k+1} =  (1-\beta)T^k + \alpha (F^k + L_{batt} \cdot \Delta u_{batt}^k)^2 + \gamma
    \label{dynaT}
\end{equation}
with $\alpha = \frac{R \cdot \Delta t}{3 V_0^2 \cdot m \cdot c  } $ and $\beta = \frac{\pi \cdot \lambda_f \cdot Nu}{m \cdot c} \Delta t$ and \[ \gamma = \frac{\Delta t}{m \cdot c}(\alpha_S \cdot D \cdot I_T + \pi \cdot \lambda_f \cdot N_u \cdot T_a + R \frac{Q_0^2}{3 V_0^2}).\] The matrix $L_{batt}$ contains the Power Transfer Distribution Factors from equation \eqref{power flow}. The signal $w^k$ represents the disturbances, as for instance the variations in power flows and temperatures due to wind evolution (variation in generation and cooling effect). It is assumed that $w^k \in W$ and $W$ is bounded. 

The aggregated state dynamics can be expressed in the following non-linear compact form:
\begin{equation}
x^{k+1} = f(x^{k},u^k,w^{k})
\end{equation}
The constraints are constituted by limitations on temperature, as well as bounds on battery capacities $E^{min,max}$ and $u_{batt}^{min,max}$, and bounds on generation curtailment. The state and control constraints can be written in the form:
\begin{equation}
H_x \cdot x +H_u \cdot u +H_w \cdot w\leq H
\label{constraintes1}
\end{equation}
Disturbances $w$ are unknown but bounded. Eq. \eqref{constraintes1} can be written to consider the worst case as:
\begin{equation}
H_x \cdot x +H_u \cdot u \leq \bar{H}
\label{constraintes}
\end{equation}
We define $P$, the mix admissible set for control and state:
\begin{equation}
P=\left\{\begin{bmatrix}
x\\
u\\
\end{bmatrix} \in \mathbb{R}^{n} \times \mathbb{R}^{n_g+1} : H_x \cdot x +H_u \cdot u \leq \bar{H}\right\}
\label{P}
\end{equation}
The goal of the controller is to stay within admissible temperatures while minimizing control costs. Only control terms are penalized in the value function: the state constraints are taken into account explicitly and doesn't represent an operational objective:
\[V_N(x,u) = \sum_{t \in \{0,N-1\}} ||\Delta u||_Q\]
with $Q$ a diagonal matrix containing the costs for curtailment and battery charge. The general MPC formulation for temperature management on electrical transmission lines using batteries and generation curtailment is:
\begin{equation}
\begin{aligned}
\label{MPCformulation1}
& V_N^0(x) = \underset{u}{\text{min}}
& & V_N(x,u) \\
& \text{s.t.}
& & x^{k+1} = f(x^{k},u^k,w^{k}), \; k \in \{0,N-1\},  \\
&&& H_x \cdot x^{k+t} + H_u \cdot u^{k} \leq \bar{H}, \; k \in \{0,N\}, \\
\end{aligned}
\end{equation}

\section{Robust problem formulation based on Tube MPC}

The formulation \eqref{MPCformulation1} contains weather data, as the Nusselt number related to the wind power and direction and the ambiant temperature. If these data are updated periodically (measurements on wind and power flows), it is possible to define several functions $f_k$ corresponding to the updates in terms of a parameter (or time) dependant model. The strategy employed here is to consider only one function $f$, and to extend the disturbance set $W$ to cope with the model mismatch. The set $W$ will handle errors related to weather variations. The simultaneous presence of constraints and disturbances can lead to infeasibility. The Tube MPC method, described in \cite{mayne2005robust}, enforces constraints satisfaction while maintaining a manageable computational complexity. The method relies on inserting suitable constraints restrictions. The tightened constraints are defined through the help of a disturbance invariant set. The method is built on linear time-invariant models of the dynamics. In order to apply this approach, a linearization of $f$ will be used next.

\subsection{Linearization of the state dynamics}

The state dynamic equation is non-linear, the nonlinearity caused by the square on the active power term $F$ anc control $\Delta u{batt}$ in \eqref{dynaT}. Let $F_0$ be the average power flow on the line. We linearize the equation \eqref{dynaT} around $F_0$ using the approximation $F =  F_0+\Tilde F \approx F_0$. The variations $\Tilde{F}$ around $F_0$ can be considered small as they are usually less than $10\%$ of the $F_0$ value.

\begin{multline}
(F +L_{batt} \Delta u_{batt})^2= (F_0 + \tilde F + L_{batt}  \Delta u_{batt})^2 \\
\approx F_0^2 + 2 F_0 \tilde F  + 2F_0 L_{batt} \Delta u_{batt} +  (L_{batt} \Delta u_{batt})^2
\\
\approx F_0^2 + 2F_0\Tilde{F}+\Tilde{u}_{batt} \\
\end{multline}
with $\Tilde{u}_{batt} =  2 F_0 L_{batt} \Delta u_{batt} +  (L_{batt} \Delta u_{batt})^2$, the new control. It can be checked that $\Tilde{u}_{batt}$ realizes a bijection with $u_{batt}$ on its domain $[-\bar{P},\bar{P}]$ if $F_0 \geq 1MW$. The assumption on $F_0$ is not strong since the usual flows are around $60-90MW$. The variations of $\Tilde{u}_{batt}$ must belong to the interval $[\bar{P}^2 - 2F_0\bar{P},\bar{P}^2+2F_0\bar{P}]$.

We denote $\Tilde{T}_{ij}$ the variations of the real temperature around $T_0$, the equilibrum point for temperature corresponding to a flow $F=F_0$.

The evolution of temperature in an electrical line can thus be written as:

\begin{equation}
\Tilde{T}^{k+1} =  (1-\beta)\Tilde{T}^k +2F_0\Tilde{\alpha} \Tilde{F}^k + \Tilde{\alpha}\Tilde{u}_{batt} + w^k
\label{eqTemp}
\end{equation}

with $\Tilde{\alpha} = \frac{R}{3 \cdot m \cdot c \cdot V_0^2} \Delta t $ and $\beta = \frac{\pi \cdot \lambda_f \cdot Nu}{m \cdot c} \Delta t$. The dynamic of the linearized system is:
\begin{equation}
    \Tilde{x}^{k+1} = A \cdot \Tilde{x}^{k+1} + B \cdot \Tilde{u}^k +w^k
    \label{sysLin}
\end{equation}
with
\begin{center}
$\Tilde{x}^k=
\begin{pmatrix}
\Tilde{F}_{ij}^k &\quad \Tilde{T}_{ij}^k & \quad \Tilde{u}_{batt}^k & \quad E_{batt}^k & \quad u_{curt,n}^k
\end{pmatrix} ^T$
$\in \mathbb{R}^{n}$
\end{center}
and \begin{center}
$A = 
\begin{pmatrix}
\mathds{1} & 0 & 0 & 0 & 0 & L \\
2F_0\cdot M_F & M_T & 0 & 0 & 0 & 0 \\
0 & 0 & 1 & 0 & 0 & 0 \\
0 & 0 & \Delta t & 1 & 0 & 0 \\
0 & 0 & 0 & 0 & 0 & \mathds{1} \\
\end{pmatrix}$,
$B = 
\begin{pmatrix}
L_{batt} & L_{curt} \\
M_F & 0 \\
1 & 0 \\
\Delta t & 0 \\
0 & 1 \\
\end{pmatrix}$
\end{center}
The matrices $M_F$ and $M_T$ are diagonal and contain the coefficients described in \eqref{eqTemp}.
\begin{center}
$M_F = 
\begin{pmatrix}
\Tilde{\alpha} & 0 & ... \\
0 & \Tilde{\alpha} & 0 \\
... & 0 & ... \\
\end{pmatrix}$
$M_T = 
\begin{pmatrix}
1-\beta & 0 & ... \\
0 & 1-\beta & 0 \\
... & 0 & ... \\
\end{pmatrix}$
\end{center}

The two matrices $L$ contain the Power Transfer Distribution Factors from equation \eqref{power flow} in the correct order. $L_{batt}$ contains the PTDF for the battery node, $L_{curt}$ contains the PTDF for the curtailment nodes.

\subsection{Robust Invariant Set Calculation}

Let $\Omega$ a positive disturbance invariant set for the system:
\begin{equation}
x^{k+1} = \bar{A}x^k+w, w \in W
\label{systemEx}
\end{equation}
$\Omega$ satisfies: \[\bar{A}x + w \in \Omega, \forall x \in \Omega, x \in W \] 

\cite{olaru2010positive} revisits the construction of robust positive invariant sets (RPI) for linear systems with additive bounded disturbances. The following theorem from \cite{kofman2007systematic} provides a simple choice for the RPI set:
\vspace{1mm}

\textbf{Theorem}:\\
Considering \eqref{systemEx}, let $\bar{A}=V^{-1}\Delta V$ be the Jordan decomposition of $\bar{A}$, a strictly stable matrix. If $W$ is described by $\bar{w}$ such that $|w| \leq \bar{w}, \forall w \in W$, then the set:
\begin{equation}
\Omega = \{x \in \mathbb{R}: |V^{-1}x|\leq |I-\Delta|^{-1}|V^{-1}B|\bar{w} + \theta\}
\end{equation}
is a RPI set with $\theta$ any arbitrary small vector with positive elements.

In the particular case of temperature management controller, we consider a closed-loop with a stabilizing gain $K$. The matrix $\bar{A}$ from the theorem represents $A+BK$ with $A$ and $B$ the matrices from \eqref{MPCformulation1}. The gain $K$ can be obtained with the pole placement method, the condition being $A+BK$ strictly stable with real eigenvalues. Obtaining a small disturbance invariant set is an important issue: the smaller the set is, the less important is the constraints tightening procedure for the tube MPC algorithm. Getting the minimum disturbance invariant set enables us to be less conservative. \cite{olaru2010positive} proposes an algorithm to refine the disturbance invariant set obtained in the theorem.

\subsection{Nominal problem and algorithm}

The robust model predictive controller from \cite{mayne2005robust} is used to guarantee the recursive feasibility for linear time-invariant prediction models and will be adapted to the temperature management problem.  The controller is based on the resolution of a nominal MPC problem. This problem is \eqref{MPCformulation1} with the tightened constraints:
\begin{equation}
\begin{bmatrix}
x\\
u\\
\end{bmatrix} \in P \ominus 
\begin{bmatrix}
1\!\!1\\
K\\
\end{bmatrix} \Omega
\label{constraintTigh}
\end{equation}

$X \ominus Y$ denotes the Pontryagin difference, defined by $X \ominus Y = \{x | x \oplus Y \in X \}$, with $X \oplus Y = \{x+y | x\in X, y \in Y \}$ the Minkowsky set addition. The robust constraint satisfaction is guaranteed if the optimization-based control problem $\mathbb{P}_N^{*}(x)$ \eqref{MPCformulationDis} is deemed feasible. The control law \eqref{law} leads to a robust closed loop functioning.
\begin{equation}
\begin{aligned}
\label{MPCformulationDis}
& V_N^{*0}(x) = \underset{u,x_0}{\text{min}}
& & V_N(x^k,u) \\
& \text{s.t.}
& & x^{k+1} = Ax^{k}+Bu^k \\
&&& \quad\quad\quad + w^{k}, \; k \in \{0,N-1\}, \\
&&& u_k \in \mathbb{U} \ominus K\Omega, \; k \in \{0,N-1\} \\
&&& x^k \in \mathbb{X} \ominus \Omega, \; k \in \{0,N-1\}\\
&&& x^k \in x_0 \oplus \Omega
\end{aligned}
\end{equation}
where $x^k$ is the state of the system at time $k$. The solution of $\mathbb{P}_N^{*}(x)$ gives the optimal control sequence $u^*(x) = u_0^*(x),...u_{N-1}^*(x)$ and the optimal state sequence $x^*(x) = x_0^*(x),...x_{N}^*(x)$. $x_0^*(x)$ is not necessary equal to $x$. The robust law applied to the system is: \begin{equation}
\kappa_N^0(x^k) = u^*_0(x^k) + K(x^k-x_0^*(x))
\label{law}
\end{equation}

The following algorithm guarantees the recursive feasibility of the MPC problem and the satisfaction of the temperature constraints.

\textbf{Algorithm}:
\begin{enumerate}
\item Define a controller $K$ such that the matrix $A_K = A+BK$ is strictly stable and calculate the associated robust invariant set $\Omega$.
\item Measure $x^k$, the power flow and temperature on the controlled lines.
\item Solve the nominal problem $\mathbb{P}_N^{*}(x^k)$ giving the nominal control $u^*_k$ and the optimal initial state $x_k^*$.
\item Apply the control $\kappa_N^k(x) = u^*_k(x) + K(x-x_k^*(x))$.
\item $k = k+1$ and go back to (2)
\end{enumerate}

\section{Simulations}

\subsection{Description of the zone}

\begin{figure}[h]
\centering
\includegraphics[scale=0.25]{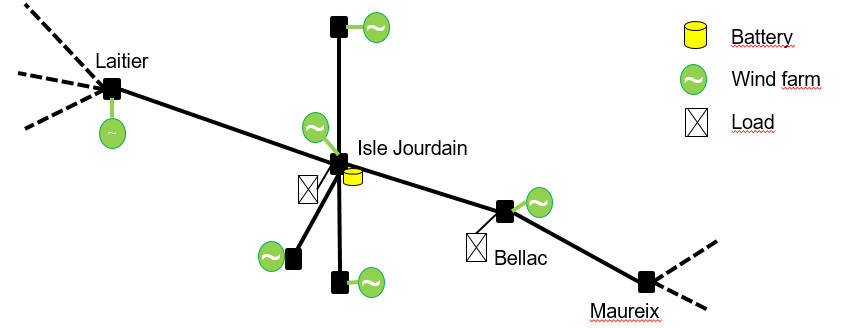}
\caption{Map of 90 kV Isle Jourdain zone}
\label{map}
\end{figure}

The robust model predictive controller for temperature management has been simulated on a zone of the French transmission network in the West area of France, see Fig. \ref{map}. Congestion are expected between Laitier and Maureix due to upcoming wind farm connections. A battery is planned to be installed in Isle Jourdain in 2020 for congestion management. \emph{Convergence}, a RTE simulation tool for network analysis, provides data to perform the simulations, see  \cite{josz2016ac} for a short presentation of this tool. We consider two lines in the zone: Isle Jourdain - Bellac and Bellac - Maureix. The conductors for both lines are Aster 288. Their characteristics are presented in \cite{refTechRTE2005} and summarized in Table \ref{characteristics}. The control is constituted by a battery and two curtailment sites: Isle Jourdain and Bellac, as their impact on the flow on the controlled lines is the biggest. Table \ref{PTDF} presents the PTDF used. It is considered that $30MW$ can be curtailed in Bellac and the battery has a $(15MW - 30MWh)$ capacity.

\begin{table}[h]
\centering
\begin{tabular}{|c|c|c|c|}
  \hline
  $m$ $(kg.m^{-1})$ & $c$ $(J.kg^{-1}.°C^{-1})$ & $D$ $(mm)$ & $R$ $(\Omega)$ \\
  \hline
  $0.627$ & $909$ & $19.6$ & $1.15e^{-4}$\\
  \hline
\end{tabular}
\caption{Aster 288 conductor characteristics}
\label{characteristics}
\end{table}
\begin{table}[h]
\centering
\begin{tabular}{|c|c|}
  \hline
  Isle Jourdain on Isle Jourdain - Bellac & 0.36\\
  \hline
  Isle Jourdain on Bellac - Maureix & 0.36\\
  \hline
  Bellac on Isle Jourdain - Bellac & 0.38\\
  \hline
  Bellac onBellac - Maureix & 0.62\\
  \hline
\end{tabular}
\caption{Power Transfer Distribution Factors}
\label{PTDF}
\end{table}

\begin{table}[h]
\centering
\begin{tabular}{|c|c|c|c|c|c|}
  \hline
  Nu & $T_a$ $(\deg C)$  & $I_T$ $(W.m^{-2})$ & $V (kV)$ & $Q(Var)$ & $F_0(MW)$\\
  \hline
  34 & 20 & 10 & 90 & 5 & 70 \\
  \hline
\end{tabular}
\caption{Weather and electrical data}
\label{disturbance}
\end{table}

Weather and electrical data considered in the simulations are presented in Table \ref{disturbance}. The method used to calculate the Nusselt number comes from \cite{iglesias2014guide}. The air thermal conductivity considered is $\lambda_f = 2.61 \cdot 10^{-2}$. The maximal disturbance for power flows is $0.1 MW$ and $0.05 \deg C$ for temperature. 

\subsection{Results}
\begin{figure}[h!]
\centering
\includegraphics[scale=0.3]{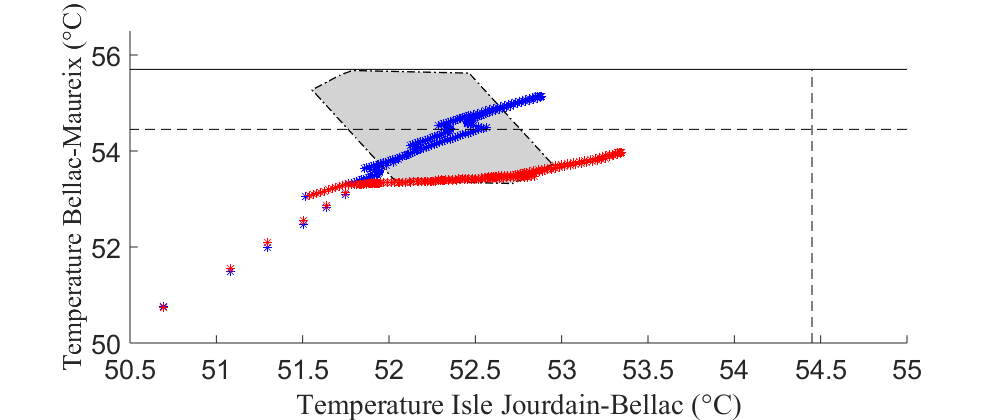}
\caption{Temperature evolution on the two controlled lines}
\label{temp2}
\end{figure}

\begin{figure}[h!]
\centering
\includegraphics[scale=0.28]{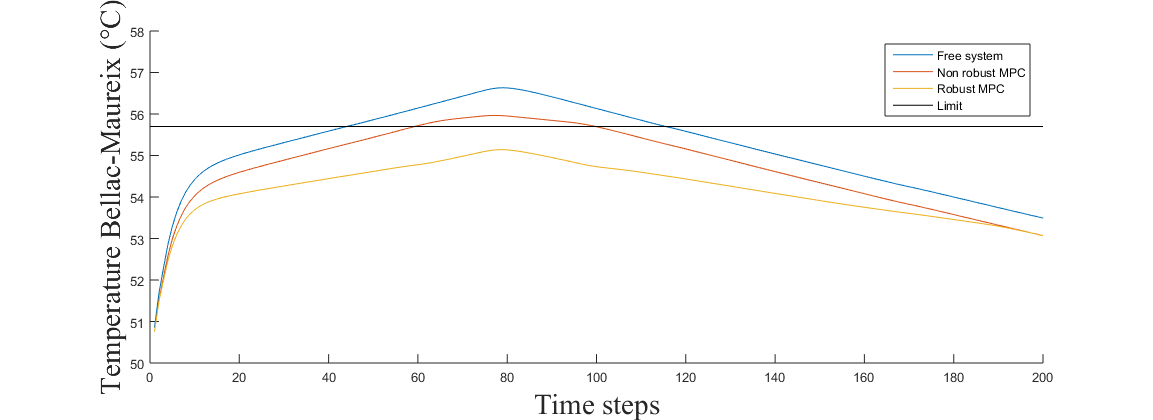}
\caption{Evolution of Bellac-Maureix line temperature}
\label{temp1}
\end{figure}

It is important to notice that the controller contains a linear model for temperature evolution, but simulations are performed with a non linear model. The point $F_0$ chosen for the Joule heating linearization is $78 MW$. The first step is to define the controller $K$ given a small disturbance invariant set $\Omega$ while maintaining the control set big enough. The tightened constraints \eqref{constraintTigh} couples constraints on control and state. For simplicity, we assume the constraints defined independently on input and state:
\[u \in \mathbb{U} \ominus K\Omega, \,  x \in \mathbb{X} \ominus \Omega\]
We consider a controller whose poles are placed at $(0.7, 0.9, 0.45, 0.21)$. Fig. \ref{temp2} shows the temperature on both lines. The set in grey is the projection of the invariant set $\Omega$ on the temperatures subspace. The maximal allowed temperature is $55.7\deg C$. The temperature on Bellac-Maureix is close to its limit, and the controller acts to reduce it. The margin required by the constraint $x_k \in \mathbb{X} \ominus \Omega$ in \eqref{MPCformulationDis} is represented by the dashed line. The horizon length is 10 prediction steps. Blue points represent the evolution of the temperature, while red points are the optimal initial state $x^*$ given by the resolution of the nominal MPC problem.
Fig. \ref{temp1} enables us to compare the temperature evolution on the free system and the controlled system. Without any corrective control action, the temperature exceeds the limit. The simulation also emphasizes the benefits of the robust controller with respect to the classical MPC one. The controller without the robust part also exceeds slightly the limit. The robust controller maintains the temperature below $55.7\deg C$, with a margin, related to the shape of the disturbance invariant set.

\begin{figure}[h!]
\centering
\includegraphics[scale=0.3]{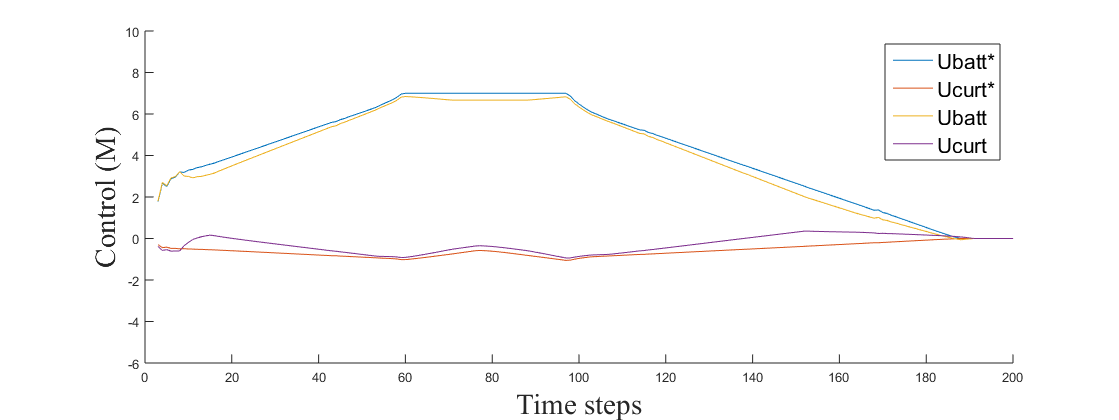}
\caption{Control evolution}
\label{control}
\end{figure}

Fig. \ref{control} shows the controls applied to the battery and generation curtailment. They are of two types: $Ubat*$ and $Ucurt*$ are the solution of \eqref{MPCformulationDis}. $Ubat^*$ and $Ucurt^*$ increases to maintain the maximal temperature below limitation in the controlled system and compensate the increase of the temperature in the free system. $Ubat$ and $Ucurt$ are the controls applied to the system taking into account the robust controller $K$ part. 

Despite the model approximation concerning the linearization of the Joule heating in the MPC formulation, the problem remains feasible and limits are not crossed, this being related to the constraints tightening approach. However, the disturbance set $W$ considered in these simulations is not truly representative of the power flows variations: $1MW$ would be a more accurate value for disturbances on power flows, especially because of errors in measurements, but considering such value results in a conservative disturbance invariant set $\Omega$, and eventually the design of a linear controller $K$ corresponding to a non-empty set $\mathbb{U} \ominus K\Omega$ becomes challenging.

\section{Conclusion}

This work presents a robust MPC controller for temperature management. The controller maximizes lines capacities as it is based on real-time lines heating, and allows therefore more renewable integration. In the model developed, the matrices $A$ and $B$ are not evolving through time, meaning that the weather conditions and the linearized point for the Joule heating are considered constant. The increase of the disturbance set is used to deal with these errors. Simulations show that the framework is robust to these modeling approximations. 

\bibliography{mybib}             
                                                   






\end{document}